# A STOCHASTIC PROCESS APPROACH TO FALSE DISCOVERY CONTROL

By Christopher Genovese[1] and Larry Wasserman[2]

## Carnegie Mellon University

This paper extends the theory of false discovery rates (FDR) pioneered by Benjamini and Hochberg [*J. Roy. Statist. Soc. Ser. B* **57** (1995) 289–300]. We develop a framework in which the False Discovery Proportion (FDP)—the number of false rejections divided by the number of rejections—is treated as a stochastic process. After obtaining the limiting distribution of the process, we demonstrate the validity of a class of procedures for controlling the False Discovery Rate (the expected FDP). We construct a confidence envelope for the whole FDP process. From these envelopes we derive confidence thresholds, for controlling the quantiles of the distribution of the FDP as well as controlling the number of false discoveries. We also investigate methods for estimating the $p$-value distribution.

## Contents

Received January 2002; revised August 2003.

[1]Supported by NSF Grant SES-98-66147.

[2]Supported by NIH Grants R01-CA54852-07 and MH57881 and NSF Grants DMS-98-03433 and DMS-01-04016.

*AMS 2000 subject classifications.* 62H15, 62G10.

*Key words and phrases.* Multiple testing, $p$-values, false discovery rate.

**Notation index**

The following summarizes the most common recurring notation and indicates where each symbol is defined.

| Symbol | Description | Section | Page |
|--------|-------------|---------|------|
| $m$ | Total number of tests performed | 2.1 | 4 |
| $P^m$ | Vector of $p$-values $(P_1, \ldots, P_m)$ | 2.2 | 4 |
| $H^m$ | Vector of hypothesis indicators $(H_1, \ldots, H_m)$ | 2.2 | 4 |
| $P_{(i)}$ | The $i$th smallest $p$-value; $P_{(0)} \equiv 0$ | 2.2 | 4 |
| $M_0$ | Number of true null hypotheses | 2.2 | 5 |
| $M_1$ | Number of false null hypotheses | 2.2 | 5 |
| $a$ | Probability of a false null | 2.2 | 4 |
| $F, f$ | Alternative $p$-value distribution (CDF, PDF) | 2.2 | 5 |
| $G, g$ | Marginal distribution (CDF, PDF) of the $P_i$'s | 2.2 | 5 |
| $\widehat{G}$ | Generic estimator of $G$ | 3 | 8 |
| $\mathbb{G}_m$ | Empirical CDF of $P^m$ | 3 | 8 |
| $U$ | Uniform CDF | 2.2 | 5 |
| $\Gamma$ | FDP process | 2.5 | 7 |
| $\Xi$ | FNP process | 2.5 | 7 |
| $\varepsilon_m$ | Dvoretzky–Kiefer–Wolfowitz nghd. radius | 3 | 8 |
| $Q$ | Asymptotic mean of $\Gamma$ | 2.5 | 7 |
| $\widetilde{Q}$ | Asymptotic mean of $\Xi$ | 2.5 | 7 |

We use $\mathbb{1}\{\ldots\}$ and $\mathsf{P}\{\ldots\}$ to denote, respectively, the indicator and probability of the event $\{\ldots\}$; subscripts on $\mathsf{P}$ specify the underlying distributions when necessary. We also use $E$ to denote expectation, and $X_m \rightsquigarrow X$ to denote that $X_m$ converges in distribution to $X$. We use $z_\alpha$ to denote the upper $\alpha$-quantile of a standard normal.

**1. Introduction.** Among the many challenges raised by the analysis of large data sets is the problem of multiple testing. In some settings it is not unusual to test thousands or even millions of hypotheses. Examples include function magnetic resonance imaging, microarray analysis in genetics and source detection in astronomy. Traditional methods that provide strong control of familywise error often have low power and can be unduly conservative in many applications.

Benjamini and Hochberg (BH) (1995, 2000) pioneered an alternative. Define the False Discovery Proportion (FDP) to be the number of false rejections divided by the number of rejections. The False Discovery Rate (FDR) is the expected FDP. BH (1995) provided a distribution-free, finite sample method for choosing a $p$-value threshold that guarantees that the FDR is less than a target level $\alpha$. The same paper demonstrated that the BH procedure is often more powerful than traditional methods that control familywise error.

Recently there has been much further work on FDR. We shall not attempt a complete review here but mention the following. Benjamini and Yekutieli (2001) extended the BH method to a class of dependent tests. Efron, Tibshirani, Storey and Tusher (2001) developed an empirical Bayes approach to multiple testing and made interesting connections with FDR. Storey (2002, 2003) connected the FDR concept with a certain Bayesian quantity and proposed a new FDR method which has higher power than the original BH method. Finner and Roters (2002) discussed the behavior of the expected number of type I errors. Sarkar (2002) considered a general class of stepwise multiple testing methods.

Genovese and Wasserman (2002) showed that, asymptotically, the BH method corresponds to a fixed threshold method that rejects all $p$-values less than a threshold $u^*$, and they characterized $u^*$. They also introduced the False Nondiscovery Rate (FNR) and found the optimal threshold $t^*$ in the sense of minimizing FNR subject to a bound on FDR. The two thresholds are related by $u^* < t^*$, implying that BH is (asymptotically) conservative. Abramovich, Benjamini, Donoho and Johnstone (2000) established a connection between FDR and minimax point estimation. (An interesting open question is whether the asymptotic results obtained in this paper can be extended to the sparse regime in the aforementioned paper where the fraction of alternatives tends to zero.)

In this paper we develop some large-sample theory for FDRs and present new methods for controlling quantiles of the false discovery distribution. An essential idea is to view the proportion of false discoveries as a stochastic process indexed by the $p$-value threshold. The problem of choosing a threshold then becomes a problem of controlling a stochastic process. Although this stochastic process is not observable, we will show that it is amenable to inference.

The main contributions of the paper include the following:

1. Development of a stochastic process framework for FDP.

|  | $H_0$ **Not Rejected** | $H_0$ **Rejected** | **Total** |
|---|---|---|---|
| $H_0$ True | $M_{0\|0}$ | $M_{1\|0}$ | $M_0$ |
| $H_0$ False | $M_{0\|1}$ | $M_{1\|1}$ | $M_1$ |
| Total | $m - R$ | $R$ | $m$ |

2. Investigation of estimators of the $p$-value distribution, even in the non-identifiable case.
3. Proof of the asymptotic validity of a class of methods for FDR control.
4. Two methods for constructing confidence envelopes for the False Discovery process and the number of false discoveries.
5. New methods, which we call *confidence thresholds*, for controlling quantiles of the false discovery distribution.

## 2. Preliminaries.

2.1. *Notation.* Consider a multiple testing situation in which $m$ tests are being performed. Suppose $M_0$ of the null hypotheses are true and $M_1 = m - M_0$ null hypotheses are false. We can categorize the $m$ tests in the following $2 \times 2$ table on whether each null hypothesis is rejected and whether each null hypothesis is true:

We define the FDP and the FNP by

$$(1) \qquad \text{FDP} = \begin{cases} \dfrac{M_{1|0}}{R}, & \text{if } R > 0, \\ 0, & \text{if } R = 0 \end{cases}$$

and

$$(2) \qquad \text{FNP} = \begin{cases} \dfrac{M_{0|1}}{m - R}, & \text{if } R < m, \\ 0, & \text{if } R = m. \end{cases}$$

The first is the proportion of rejections that are incorrect, and the second—the dual quantity—is the proportion of nonrejections that are incorrect. Notice that $\text{FDR} = \mathsf{E}(\text{FDP})$, and following Genovese and Wasserman ([2002](#)), we define $\text{FNR} = \mathsf{E}(\text{FNP})$. Storey ([2002](#)) considered a different definition of FDR, called pFDR for positive FDR, by conditioning on the event that $R > 0$ and discussed the advantages and disadvantages of this definition.

2.2. *Model.* Let $H_i = 0$ (or 1) if the $i$th null hypothesis is true (false) and let $P_i$ denote the $i$th $p$-value. Define vectors $P^m = (P_1, \ldots, P_m)$ and $H^m = (H_1, \ldots, H_m)$. Let $P_{(1)} < \cdots < P_{(m)}$ denote the ordered $p$-values, and define $P_{(0)} \equiv 0$.

In this paper we use a random effects (or hierarchical) model as in Efron, Tibshirani, Storey and Tusher ([2001](#)). Specifically, we assume the following for $0 \leq a \leq 1$:

$$H_1, \ldots, H_m \sim \text{Bernoulli}(a),$$

$$\Xi_1, \ldots, \Xi_m \sim \mathcal{L}_{\mathcal{F}},$$

$$P_i | H_i = 0, \quad \Xi_i = \xi_i \sim \text{Uniform}(0, 1),$$

$$P_i | H_i = 1, \quad \Xi_i = \xi_i \sim \xi_i,$$

where $\Xi_1, \ldots, \Xi_m$ denote distribution functions and $\mathcal{L}_{\mathcal{F}}$ is an arbitrary probability measure over a class of distribution functions $\mathcal{F}$ that stochastically dominates the Uniform$(0,1)$.

It follows that the marginal distribution of the $p$-values is

$$(3) \qquad\qquad G = (1-a)U + aF,$$

where $U(t)$ denotes the Uniform$(0,1)$ CDF and $F(t) = \int \xi(t) \, d\mathcal{L}_{\mathcal{F}}(\xi)$. Note that $G \geq U$. Except where noted we assume that $G$ is strictly concave with density $g = G'$.

REMARK 2.1. A more common approach in multiple testing is to use a conditional model in which $H_1, \ldots, H_m$ are fixed, unknown binary values. The results in this paper can be cast in a conditional framework but we find the random effects framework to be more intuitive.

Define $M_0 = \sum_i (1 - H_i)$ and $M_1 = \sum_i H_i$. Hence, $M_0 \sim \text{Binomial}(m, 1-a)$ and $M_1 = m - M_0$.

2.3. *The Benjamini–Hochberg and plug-in methods.* The Benjamini–Hochberg (BH) procedure is a distribution-free method for choosing which null hypotheses to reject while guaranteeing that FDR $\leq \alpha$ for some preselected level $\alpha$. The procedure rejects all null hypotheses for which $P_i \leq P_{(R_{\mathrm{BH}})}$, where

$$(4) \qquad\qquad R_{\mathrm{BH}} = \max\left\{ 0 \leq i \leq m : P_{(i)} \leq \alpha \frac{i}{m} \right\}.$$

BH ([1995](#)) proved that this procedure guarantees

$$(5) \qquad\qquad \mathsf{E}(\text{FDP} \,|\, M_0) \leq \frac{M_0}{m}\alpha \leq \alpha,$$

regardless of how many nulls are true and regardless of the distribution of the $p$-values under the alternatives. (When the $p$-value distribution is continuous, BH shows that the first inequality is an equality.) In the context of our model, this result becomes

$$(6) \qquad\qquad \text{FDR} \leq (1-a)\alpha \leq \alpha.$$

Genovese and Wasserman ([2002](#)) showed that, asymptotically, the BH procedure corresponds to rejecting the null when the $p$-value is less than $u^*$, where $u^*$ is the solution to the equation $G(u) = u/\alpha$, in the notation of the current paper. This $u^*$ satisfies $\alpha/m \leq u^* \leq \alpha$ for large $m$, which shows that the BH method is intermediate between Bonferroni (corresponding to $\alpha/m$) and uncorrected testing (corresponding to $\alpha$). They also showed that $u^*$ is strictly less than the optimal $p$-value cutoff.

Benjamini and Hochberg ([2000](#)), in work originally written in 1994, showed that the power of the BH ([1995](#)) procedure could be improved by estimating the number of true null hypotheses. They also proposed an estimator of $\widehat{\mathrm{FDR}}(t)$ and proposed a threshold $T$ that maximizes the number of rejections subject to $\widehat{\mathrm{FDR}}(T) \leq \alpha$. A similar approach was investigated in Storey ([2002](#)) and Storey, Taylor and Siegmund ([2004](#)). It remains an open question whether $\mathrm{FDR}(T) \leq \alpha$. We address an asymptotic version of this question in Section 5.

The threshold $T$ chosen this way can also be viewed as a plug-in estimator. Let

$$(7) \qquad t(a, G) = \sup \left\{ t : \frac{(1-a)t}{G(t)} \leq \alpha \right\}.$$

Suppose we reject whenever the $p$-value is less than $t(a, G)$. From Genovese and Wasserman ([2002](#)) it follows that, asymptotically, the FDR is less than $\alpha$. The intuition for (7) is that $(1-a)t/G(t)$ is, up to an exponentially small term, the FDR at a fixed threshold $t$. Moreover, if $G$ is concave this threshold has the smallest asymptotic FNR among all procedures with FDR less than or equal to $\alpha$ [cf. Genovese and Wasserman ([2002](#))]. We call $t(a, G)$ the *oracle threshold*. The standard plug-in method is to estimate the functional $t(a, G)$ by $T = t(\hat{a}, \widehat{G})$, where $\hat{a}$ and $\widehat{G}$ are estimators of $a$ and $G$. Let $\mathbb{G}_m$ be the empirical CDF of $P^m$. Theorem 2 of BH ([1995](#)) shows that $T_{\mathrm{BH}} = t(0, \mathbb{G}_m)$ yields the BH threshold. Benjamini and Hochberg ([2000](#)) and Storey ([2002](#)) showed that $T = t(\hat{a}_0, \mathbb{G}_m)$ has higher power than the BH threshold, where

$$\hat{a}_0 = \max \left( 0, \frac{\mathbb{G}_m(t_0) - t_0}{1 - t_0} \right)$$

and $t_0 \in (0, 1)$. Clearly, other estimators of $a$ and $G$ are possible and we shall call any threshold of the form $T = t(\hat{a}, \widehat{G})$ a plug-in threshold.

We describe alternative estimators of $a$ in Section [3.2](#). Storey ([2002](#)) provided simulations to show that the plug-in procedure has good power but did not provide a proof that it controls FDR at level $\alpha$. We settle this question in Section [5](#) where we show that, under weak conditions on $\hat{a}$, the procedure asymptotically controls FDR at level $\alpha$.

**2.4. *Multiple testing procedures*.**   A *multiple testing procedure* $T$ is a mapping taking $[0, 1]^m$ into $[0, 1]$, where it is understood that the null hypotheses corresponding to all $p$-values less than $T(P^m)$ are rejected. We often call $T$ the *threshold*.

Let $\alpha, t \in [0, 1]$ and $0 \leq r \leq m$, and recall that $P_{(0)} \equiv 0$. Let $\widehat{G}$ and $\hat{g}$ be generic estimates of $G$ and $g = G'$, respectively. Similarly, let $\widehat{\mathsf{P}}\{H = h | P = t\}$ denote an estimator of $\mathsf{P}\{H = h | P = t\}$.

Some examples of multiple testing procedures will illustrate the generality of the framework:

Uncorrected testing $\qquad T_U(P^m) = \alpha$

Bonferroni $\qquad T_B(P^m) = \alpha/m$

Fixed threshold at $t$ $\qquad T_t(P^m) = t$

Benjamini–Hochberg $\qquad T_{BH}(P^m) = \sup\{t : \mathbb{G}_m(t) = t/\alpha\} = P_{(R_{BH})}$

Oracle $\qquad T_o(P^m) = \sup\{t : G(t) = (1-a)t/\alpha\}$

Plug in $\qquad T_{PI}(P^m) = \sup\{t : \widehat{G}(t) = (1-a)t/\hat{\alpha}\}$

First $r$ $\qquad T_{(r)} = P_{(r)}$

Bayes' classifier $\qquad T_{BC}(P^m) = \sup\{t : \hat{g}(t) > 1\}$

Regression classifier $\qquad T_{Reg}(P^m) = \sup\{t : \widehat{\mathsf{P}}\{H_1 = 1 | P_1 = t\} > 1/2\}.$

**2.5. *FDP and FNP as stochastic processes.*** An important idea that we use throughout the paper is that the FDP, regarded as a function of the threshold $t$, is a stochastic process. This observation is crucial for studying the properties of procedures.

Define the *FDP process*

$$\text{(8)} \qquad \Gamma(t) \equiv \Gamma(t, P^m, H^m) = \frac{\sum_i \mathbb{1}\{P_i \leq t\}(1 - H_i)}{\sum_i \mathbb{1}\{P_i \leq t\} + \mathbb{1}\{\text{all} P_i > t\}},$$

where the last term in the denominator makes $\Gamma = 0$ when no $p$-values are below the threshold. Also define the *FNP process*

$$\text{(9)} \qquad \Xi(t) \equiv \Xi(t, P^m, H^m) = \frac{\sum_i \mathbb{1}\{P_i > t\} H_i}{\sum_i \mathbb{1}\{P_i > t\} + \mathbb{1}\{\text{all} P_i \leq t\}}.$$

The FDP and FNP of a procedure $T$ are $\Gamma(T) \equiv \Gamma(T(P^m), P^m, H^m)$ and $\Xi(T) \equiv \Xi(T(P^m), P^m, H^m)$. Let

$$\text{(10)} \qquad Q(t) = (1-a)\frac{t}{G(t)},$$

$$\text{(11)} \qquad \tilde{Q}(t) = a\frac{1 - F(t)}{1 - G(t)}.$$

The following lemma is a corollary of Theorem 1 in Storey (2002).

LEMMA 2.1. *Under the mixture model, for $t > 0$,*

$$\mathsf{E}\Gamma(t) = Q(t)(1 - (1 - G(t))^m),$$

$$\mathsf{E}\Xi(t) = \tilde{Q}(t)(1 - G(t)^m).$$

*The second terms on the right-hand side of both equations differ from 1 by an exponentially small quantity.*

One of the essential difficulties in studying a procedure $T$ is that $\Gamma(T)$ is the evaluation of the stochastic process $\Gamma(\cdot)$ at a random variable $T$. Both depend on the observed data, and in general they are correlated. In particular, if $\widehat{Q}(t)$ estimates $\mathrm{FDR}(t)$ well at each fixed $t$, it does not follow that $\widehat{Q}(T)$ estimates $\mathrm{FDR}(T)$ well at a random $T$. The stochastic process point of view provides a suitable framework for addressing this problem.

**3. Estimating the $p$-value distribution.** Recall that, under the mixture model, $P_1, \ldots, P_m$ have CDF $G(t) = (1-a)\,t + a\,F(t)$. Let $\widehat{G}$ denote an estimator of $G$. Let $\mathbb{G}_m$ denote the empirical CDF of $P^m$. We will use the Dvoretzky–Kiefer–Wolfowitz (DKW) inequality: for any $x > 0$,

$$\tag{12} \mathsf{P}\{\|\mathbb{G}_m(t) - G(t)\|_\infty > x\} \le 2e^{-2mx^2},$$

where $\|F - G\|_\infty = \sup_{0 \le t \le 1} |F(t) - G(t)|$. Given $\alpha \in (0, 1)$, let

$$\tag{13} \varepsilon_m \equiv \varepsilon_m(\alpha) = \sqrt{\frac{1}{2m} \log\left(\frac{2}{\alpha}\right)}$$

so that, from DKW, $\mathsf{P}\{\|\mathbb{G}_m(t) - G(t)\|_\infty > \varepsilon_m\} \le \alpha$.

Several improvements on $\mathbb{G}_m$ are possible. Since $G \ge U$, we replace any estimator $\mathbb{G}_m$ with $\max\{\mathbb{G}_m(t), t\}$. When $G$ is assumed to be concave, a better estimate of $G$ is the least concave majorant (LCM) $\mathbb{G}_{\mathrm{LCM},\,m}$ defined to be the infimum of the set of all concave CDF's lying above $\mathbb{G}_m$. Most $p$-value densities in practical problems are decreasing in $p$, which implies that $G$ is concave. We can also replace $\mathbb{G}_{\mathrm{LCM},m}$ with $\max\{\mathbb{G}_{\mathrm{LCM},m}(t), t\}$. The DKW inequality and the standard limiting results still hold for the modified versions of both estimators. We will thus use $\widehat{G}$ to denote the modified estimators in either case. We will indicate explicitly if concavity is required or if the LCM estimator is proscribed.

Once we obtain estimates $\hat{a}$ and $\widehat{G}$, we define

$$\tag{14} \widehat{Q}(t) = \frac{(1 - \hat{a})}{\widehat{G}(t)}.$$

3.1. *Identifiability and purity.* Before discussing the estimation of $a$, it is helpful to first discuss identifiability. For example, if $a$ is not identifiable, there is no guarantee that the estimate used in the plug-in method will give good performance. The results in the ensuing sections show that despite not being completely identified, it is possible to make sensible inferences about $a$.

Say that $F$ is *pure* if $\mathrm{ess\,inf}_t\, f(t) = 0$, where $f$ is the density of $F$. Let $\mathcal{O}_F$ be the set of pairs $(b, H)$ such that $b \in [0, 1]$, $H \in \mathcal{F}$ and $F = (1 - b)U + bH$. $F$ is identifiable if $\mathcal{O}_F = \{(1, F)\}$.

Define

$$\zeta_F = \inf\{b : (b, H) \in \mathcal{O}_F\},$$

$$\underline{F} = \frac{F - (1 - \zeta_F)U}{\zeta_F},$$

$$\underline{a}_F = a\zeta_F.$$

We will often drop the subscript $F$ on $\underline{a}_F$ and $\zeta_F$. Note that $G$ can be decomposed as

$$G = (1 - a)U + aF$$

$$= (1 - a)U + a[(1 - \zeta)U + \zeta \underline{F}]$$

$$= (1 - a\zeta)U + a\zeta \underline{F}$$

$$= (1 - \underline{a})U + \underline{a}\,\underline{F}.$$

Purity implies identifiability but not vice versa. Consider the following example. Let $\mathcal{F}$ be the Normal $(\theta, 1)$ family and consider testing $H_0 : \theta = 0$ versus $H_1 : \theta \neq 0$. The density of the $p$-value is

$$f_\theta(p) = \tfrac{1}{2}e^{-n\theta^2/2}[e^{-\sqrt{n}\theta\Phi^{-1}(1 - p/2)} + e^{\sqrt{n}\theta\Phi^{-1}(1 - p/2)}].$$

Now, $f_\theta(1) = e^{-n\theta^2/2} > 0$ so this test is impure. However, the parametric assumption makes $a$ and $\theta$ identifiable when the null is false. It is worth noting that $f_\theta(1)$ is exponentially small in $n$. Hence, the difference between $a$ and $\underline{a}$ is small. Even when $X$ has a $t$-distribution with $\nu$ degrees of freedom, $f_\theta(1) = (1 + n\theta^2/\nu)^{-(\nu+1)/2}$. Thus, in practical cases, $a - \underline{a}$ will be quite small. On the other hand, one-sided tests for continuous exponential families are pure and identifiable.

The problem of estimating $a$ has been considered by Efron, Tibshirani, Storey and Tusher (2001) and Storey (2002) who also discussed the identifiability issue. In particular, Storey noted that $G(t) = (1 - a)t + aF(t) \leq (1 - a)t + a$ for all $t$. It then follows that, for any $t_0 \in (0, 1)$,

$$(15) \qquad 0 \leq a_0 \equiv \frac{G(t_0) - t_0}{1 - t_0} \leq \underline{a} \leq a \leq 1.$$

Thus, an identifiable lower bound on $a$ is $a_0$. The following result gives precise information about the best bounds that are possible.

PROPOSITION 3.1. *If $F$ is absolutely continuous and stochastically dominates $U$, then*

$$\zeta = 1 - \inf_t F'(t) \quad and \quad \underline{a} = 1 - \inf_t G'(t).$$

*If $F$ is concave, then the infima are achieved at $t = 1$. For any $b \in [\zeta, 1]$ we can write $G = (1 - ab)U + abF_b$, where $F_b = (F - (1 - b)U)/b$ is a CDF and $F \leq F_b$.*

3.2. *Estimating $a$.* Here we discuss estimating $a$. Related work includes Schweder and Spjøtvoll (1982), Hochberg and Benjamini (1990), Benjamini and Hochberg (2000) and Storey (2002).

We begin with a uniform confidence interval for $\underline{a}$.

THEOREM 3.1. *Let*

$$(16) \qquad a_* = \max_t \frac{\widehat{G}(t) - t - \varepsilon_m}{1 - t}.$$

*Then $[a_*, 1]$ is a uniform $1 - \alpha$ confidence interval for $\underline{a}$, that is,*

$$(17) \qquad \inf_{a, F} \mathsf{P}_{a,F}\{\underline{a} \in [a_*, 1]\} \geq 1 - \alpha,$$

*and if one restricts $\widehat{G}$ to be the empirical distribution function, then for each $(a, F)$ pair,*

$$(18) \qquad \mathsf{P}_{a,F}\{\underline{a} \in [a_*, 1]\} \leq 1 - \alpha + 2 \sum_{j=1}^{\infty} (-1)^{j+1} \left(\frac{\alpha}{2}\right)^{j^2} + O\left(\frac{(\log m)^2}{\sqrt{m}}\right),$$

*where the remainder term may depend on $a$ and $F$. Because $a \geq \underline{a}$, $[a_*, 1]$ is a valid finite-sample $1 - \alpha$ confidence interval for $a$ as well.*

PROOF. The inequality (17) follows immediately from DKW because $G(t) \geq \widehat{G}(t) - \varepsilon_m$ for all $t$ with probability at least $1 - \alpha$. The sum on the right-hand side of (18) follows from the closed-form limiting distribution of the Kolmogorov–Smirnov statistic, and the order of the error follows from the Hungarian embedding. To see this, note that

$$\underline{a} < a_* \implies \underline{a}\sqrt{m} < \max_t \sqrt{m}\frac{\mathbb{G}_m(t) - G(t)}{1 - t} + \sqrt{m}\frac{G(t) - t}{1 - t} - \frac{\varepsilon_m \sqrt{m}}{1 - t}$$

$$\implies \underline{a}\sqrt{m} < \max_t \sqrt{m}\frac{\mathbb{G}_m(t) - G(t)}{1 - t} + \sqrt{m}\underline{a} - \frac{\varepsilon_m \sqrt{m}}{1 - t}$$

$$\implies 0 < \max_t \sqrt{m}\frac{\mathbb{G}_m(t) - G(t)}{1 - t} - \frac{\varepsilon_m \sqrt{m}}{1 - t}$$

$$\implies 0 < \max_t \sqrt{m}(\mathbb{G}_m(t) - G(t)) - \varepsilon_m \sqrt{m}$$

$$\implies \|\sqrt{m}(\mathbb{G}_m(t) - G(t))\|_\infty > \varepsilon_m \sqrt{m}.$$

Hence,

$$(19) \qquad \mathsf{P}\{\underline{a} < a_*\} \leq \mathsf{P}\{\|\sqrt{m}(\mathbb{G}_m(t) - G(t))\|_\infty > \varepsilon_m \sqrt{m}\}.$$

Next apply the Hungarian embedding [van der Vaart (1998), page 269]:

$$\limsup_{m \to \infty} \frac{\sqrt{m}}{(\log m)^2} \|\sqrt{m}(\mathbb{G}_m - G) - \mathbb{B}_m\|_\infty < \infty \qquad \text{a.s.,}$$

for a sequence of Brownian bridges $\mathbb{B}_m$. Recall the distribution of the Kolmogorov–Smirnov statistic:

$$\mathsf{P}\{\|\mathbb{B}\|_\infty > x\} = 2\sum_{j=1}^{\infty}(-1)^{j+1}e^{-2j^2x^2},$$

for a generic Brownian bridge $\mathbb{B}$. The result follows by taking $x = \sqrt{m}\,\varepsilon_m$. In the concave case, the LCM can be substituted for $\widehat{G}$ and the result still holds since, by Marshall's lemma, $\|\widehat{G}_{\mathrm{LCM},m} - G\|_\infty \leq \|\widehat{G}_m - G\|_\infty$. $\quad\square$

PROPOSITION 3.2 (Storey's estimator). *Fix $t_0 \in (0,1)$ and let*

$$\hat{a}_0 = \left(\frac{\mathbb{G}_m(t_0) - t_0}{1 - t_0}\right)_+.$$

*If $G(t_0) > t_0$,*

$$\hat{a}_0 \xrightarrow{P} \frac{G(t_0) - t_0}{1 - t_0} \equiv a_0 \leq \underline{a},$$

*and*

$$\sqrt{m}\left(\hat{a}_0 - \frac{G(t_0) - t_0}{1 - t_0}\right) \rightsquigarrow N\left(0, \frac{G(t_0)(1 - G(t_0))}{(1 - t_0)^2}\right).$$

*If $G(t_0) = t_0$,*

$$\sqrt{m}\hat{a}_0 \rightsquigarrow \frac{1}{2}\delta_0 + \frac{1}{2}N^+\left(0, \frac{t_0}{1 - t_0}\right),$$

*where $\delta_0$ is a point-mass at zero and $N^+$ is a positive-truncated normal.*

A consistent estimate of $\underline{a}$ is available if we assume weak smoothness conditions on $g$. For example, one can use the spacings estimator of Swanepoel (1999) which is of the form $2r_m/(mV_m)$, where $r_m = m^{4/5}(\log m)^{-2\delta}$ and $V_m$ is a selected spacing in the order statistics of the $p$-values.

THEOREM 3.2. *Assume that at the value $t$ where $g$ achieves its minimum, $g''$ is bounded away from 0 and $\infty$ and Lipschitz of order $\lambda > 0$. For every $\delta > 0$, there exists an estimator $\hat{a}$ such that*

$$\frac{m^{(2/5)}}{(\log m)^\delta}(\hat{a} - \underline{a}) \rightsquigarrow N(0, (1 - \underline{a})^2).$$

PROOF. Let $\hat{a}$ be the estimator defined in Swanepoel (1999) with $r_m = m^{4/5}(\log m)^{-2\delta}$ and $s_m = m^{4/5}(\log m)^{4\delta}$. The result follows from Swanepoel [(1999), Theorem 1.3]. $\quad\square$

REMARK 3.1. An alternative estimator is $\hat{a} = 1 - \min_t \hat{g}(t)$, where $\hat{g}$ is a kernel estimator.

Now suppose we assume only that $G$ is concave and hence $g = G'$ is decreasing. Hengartner and Stark (1995) derived a finite-sample confidence envelope $[\gamma^-(\cdot), \gamma^+(\cdot)]$ for a density $g$ assuming only that it is monotone. Define

$$\hat{a}_{\mathrm{HS}} = 1 - \min\{h(1) : \gamma^- \le h \le \gamma^+\}.$$

THEOREM 3.3. If $G$ is concave and $g = G'$ is Lipschitz of order 1 in a neighborhood of 1, then

$$\left(\frac{n}{\log n}\right)^{1/3} (\hat{a}_{\mathrm{HS}} - \underline{a}) \xrightarrow{P} 0.$$

Also, $[1 - \gamma^+(1), 1 - \gamma^-(1)]$ is a $1 - \alpha$ confidence interval for $\underline{a}$ for $0 \le \alpha \le 1$ and all $m$. Further,

$$\inf_{a, F} \mathsf{P}\{a \in [1 - \gamma^+(1), 1]\} \ge 1 - \alpha,$$

where the infimum is over all concave $F$'s.

PROOF. Follows from Hengartner and Stark (1995). □

3.3. *Estimating $F$.* It may be useful in some cases to estimate the alternative mixture distribution $F$. There are many possible methods; we consider here projection estimators defined by

$$(20) \qquad \widehat{F}_m = \arg\min_{H \in \mathcal{F}} \|\widehat{G} - (1 - \hat{a})U - \hat{a}H\|_\infty,$$

where $\hat{a}$ is an estimate of $a$. The Appendix gives an algorithm to find $\widehat{F}_m$.

It is helpful to consider first the case where $a$ is known, and here we substitute $a$ for $\hat{a}$ in the definition of $\widehat{F}_m$.

THEOREM 3.4. Let

$$\widehat{F}_m = \arg\min_{H \in \mathcal{F}} \|\widehat{G} - (1 - a)U - aH\|_\infty.$$

Then

$$\|F - \widehat{F}_m\|_\infty \le \frac{2\|G - \widehat{G}\|_\infty}{a} \xrightarrow{\text{a.s.}} 0.$$

PROOF.

$$a\|F - \widehat{F}_m\|_\infty = \|aF - a\widehat{F}_m\|_\infty$$
$$= \|(1-a)U + aF - (1-a)U - a\widehat{F}_m\|_\infty$$
$$= \|G - (1-a)U - a\widehat{F}_m\|_\infty$$
$$= \|G - \widehat{G} + \widehat{G} - (1-a)U - a\widehat{F}_m\|_\infty$$
$$\leq \|\widehat{G} - G\|_\infty + \|\widehat{G} - (1-a)U - a\widehat{F}_m\|_\infty$$
$$\leq \|\widehat{G} - G\|_\infty + \|\widehat{G} - (1-a)U - aF\|_\infty$$
$$= \|\widehat{G} - G\|_\infty + \|\widehat{G} - G\|_\infty.$$

The last statement follows from the uniform consistency of $\widehat{G}$. $\square$

When $a$ is unknown, the projection estimator $\widehat{F}$ is consistent whenever we have a consistent estimator of $\underline{a}$. Recall that in the identifiable case $a = \underline{a}$ and $F = \underline{F}$.

THEOREM 3.5. *Let $\hat{a}$ be a consistent estimator of $\underline{a}$. Then*

$$\|\widehat{F}_m - \underline{F}\|_\infty \leq \frac{\|\widehat{G} - G\|_\infty + |\hat{a} - \underline{a}|}{\underline{a}} \xrightarrow{P} 0.$$

PROOF. Let $\delta_m = \|\widehat{G} - (1-\hat{a})U - \hat{a}\widehat{F}\|_\infty$. Since $\widehat{F}$ is the minimizer,

$$\delta_m \leq \|\widehat{G} - (1-\hat{a})U - \hat{a}\underline{F}\|_\infty$$
$$= \|\widehat{G} - G + (\hat{a} - \underline{a})U - (\hat{a} - \underline{a})\underline{F}\|_\infty$$
$$\leq \|\widehat{G} - G\|_\infty + |\hat{a} - \underline{a}|$$
$$\xrightarrow{P} 0.$$

We also have that

$$\delta_m \geq |\|\widehat{G} - (1-\hat{a})U - \hat{a}\underline{F}\|_\infty - \hat{a}\|\underline{F} - \widehat{F}\|_\infty|.$$

Since $\delta_m$ and $\|\widehat{G} - (1-\hat{a})U - \hat{a}\underline{F}\|_\infty \xrightarrow{P} 0$ by the above and $\hat{a} \xrightarrow{P} \underline{a}$, it follows that $\|\underline{F} - \widehat{F}\|_\infty \xrightarrow{P} 0$. Moreover,

$$\|\underline{F} - \widehat{F}\|_\infty \leq \frac{\|\widehat{G} - G\|_\infty + |\hat{a} - \underline{a}|}{\underline{a}}. \qquad \square$$

**4. Limiting distributions.** In this section we discuss the limiting distribution of $\Gamma$ and $\widehat{Q}$. Let

$$\Lambda_0(t) = \frac{1}{m} \sum_{i=1}^m (1 - H_i) \mathbb{1}\{P_i \le t\} \quad \text{and} \quad \Lambda_1(t) = \frac{1}{m} \sum_{i=1}^m H_i \mathbb{1}\{P_i \le t\},$$

and, for each $c \in (0, 1)$, define

$$\Omega_c(t) = (1 - c)\Lambda_0(t) - c\Lambda_1(t) = \frac{1}{m} \sum_i D_i(t),$$

where $D_i(t) = \mathbb{1}\{P_i \le t\}(1 - H_i - c)$. Let

$$\mu_c(t) = \mathsf{E} D_1(t) = (1 - a)t - cG(t).$$

Let $(W_0, W_1)$ be a continuous two-dimensional mean zero Gaussian process with covariance kernel $R_{ij}(s, t) = \mathrm{Cov}(W_i(s), W_j(t))$ given by

$$(21) \quad R(s, t) = \left[ \begin{pmatrix} (1-a)(s \wedge t) - (1-a)^2 st & -(1-a)s\,aF(t) \\ -(1-a)t\,aF(s) & aF(s \wedge t) - a^2 F(s)F(t) \end{pmatrix} \right].$$

THEOREM 4.1. *Let $W$ be a continuous mean zero Gaussian process with covariance*

$$(22) \quad \begin{aligned} K_\Omega(s, t) &= (1-a)(1-c)[(1-c)(s \wedge t - (1-a)st) + ac(tF(s) + sF(t))] \\ &\quad + ac[cF(s \wedge t) - acF(s)F(t)]. \end{aligned}$$

*Then*

$$\sqrt{m}(\Omega_c - \mu_c) \rightsquigarrow W.$$

PROOF. Let

$$Z_m(t) = \sqrt{m}(\Omega_c(t) - \mu_c(t)) \quad \text{and} \quad Z_m^*(t) = \sqrt{m}(\Omega_c^*(t) - \hat{\mu}_c(t))$$

for $t \in [0, 1]$. Let

$$(W_{m,0}(t), W_{m,1}(t)) \equiv (\sqrt{m}(\Lambda_0(t) - (1-a)t), \sqrt{m}(\Lambda_1(t) - aF(t))).$$

By standard empirical process theory, $(W_{m,0}(t), W_{m,1}(t))$ converges to $(W_0, W_1)$. The covariance kernel $R$ stated in (21) follows by direct calculation. The result for $\Omega_c$ is immediate since $\Omega_c$ is a linear combination of $\Lambda_0$ and $\Lambda_1$. $\square$

THEOREM 4.2 (Limiting distribution of FDP process). *For $t \in [\delta, 1]$ for any $\delta > 0$, let*

$$Z_m(t) = \sqrt{m}(\Gamma_m(t) - Q(t)).$$

*Let $Z$ be a Gaussian process on $(0,1]$ with mean 0 and covariance kernel*

$$K_\Gamma(s,t) = a(1-a)\frac{(1-a)stF(s \wedge t) + aF(s)F(t)(s \wedge t)}{G^2(s)G^2(t)}.$$

*Then $Z_m \rightsquigarrow Z$ on $[\delta,1]$.*

REMARK 4.1. The reason for restricting the theorem to $[\delta,1]$ is that the variance of the process is infinite at zero.

PROOF OF THEOREM 4.1. Note that $\Gamma_m(t) = \Lambda_0(t)/(\Lambda_0(t) + \Lambda_1(t)) \equiv r(\Lambda_0, \Lambda_1)$, where $\Lambda_0$ and $\Lambda_1$ are defined as before and $r(\cdot, \cdot)$ maps $\ell^\infty \times \ell^\infty \to \ell^\infty$, where $\ell^\infty$ is the set of bounded functions on $(\delta,1]$ endowed with the sup norm. Note that $r((1-a)U, aF) = Q$. It can be verified that $r(\cdot, \cdot)$ is Fréchet differentiable at $((1-a)U, aF)$ with derivative

$$r'_{((1-a)U,aF)}(V) = \frac{aFV_0 - (1-a)UV_1}{G^2},$$

where $U(t) = t$, $V = (V_0, V_1)$. Hence, by the functional delta method [van der Vaart (1998), Theorem 20.8],

$$Z_m \rightsquigarrow r'_{((1-a)U,aF)}(W) = \frac{aFW_0 - (1-a)UW_1}{G^2},$$

where $(W_0, W_1)$ is the process defined just before (21). The covariance kernel of the latter expression is $K_\Gamma(s,t)$. □

REMARK 4.2. A Gaussian limiting process can be obtained for FNP [i.e., $\Xi(t)$] along similar lines.

The next theorems follow from the previous results followed by an application of the functional delta method.

THEOREM 4.3. *Let $\widehat{Q}(t) = (1-a)t/\widehat{G}(t)$. For any $\delta > 0$,*

$$\sqrt{m}(\widehat{Q}(t) - Q(t)) \rightsquigarrow W$$

*on $[\delta,1]$, where $W$ is a mean zero Gaussian process on $(0,1]$ with covariance kernel*

$$K_Q(s,t) = Q(s)Q(t)\frac{G(s \wedge t) - G(s)G(t)}{G(s)G(t)}.$$

THEOREM 4.4. *Let $\widehat{Q}(t) = (1-a)t/\widehat{G}(t)$. We have*

$$\sqrt{m}(\widehat{Q}^{-1}(v) - Q^{-1}(v)) \rightsquigarrow W,$$

*where $W$ is a mean zero Gaussian process with covariance kernel*

$$K_{Q^{-1}}(u,v) = \frac{K_Q(s,t)}{Q'(s)Q'(t)}$$

$$= (1-a)^2 uv \frac{G(s \wedge t) - G(s)G(t)}{[1-a-ug(s)][1-a-vg(t)]},$$

*with $s = Q^{-1}(u)$ and $t = Q^{-1}(v)$.*

THEOREM 4.5.  *Let $\widehat{Q}(t) = (1 - \hat{a}_0)t/\widehat{G}(t)$, where $\hat{a}_0$ is Storey's estimator. Then*

$$\sqrt{m}(\widehat{Q}(t) - Q(t)) \rightsquigarrow W,$$

*where $W$ is a mean zero Gaussian process with covariance kernel*

$$K(s,t) = \frac{t^2}{(1-t_0)^2 G^2(s) G^2(t)}$$

$$\times \Big( G(s)G(t)t_0(1-t_0) + G(t)(1-G(t_0))R(s,t_0)$$

$$+ G(s)(1-G(t_0))R(t,t_0) + (1-G(t_0))^2 R(s,t) \Big),$$

*where $R(s,t) = s \wedge t - st$.*

**5. Asymptotic validity of plug-in procedures.**  Let $\widehat{Q}^{-1}(c) = \sup\{0 \le t \le 1 : \widehat{Q}(t) \le c\}$. Then the plug-in threshold $T_{\text{PI}}$ defined earlier can be written $T_{\text{PI}}(P^m) = \widehat{Q}^{-1}(\alpha)$. Here we establish the asymptotic validity of $T_{\text{PI}}$ in the sense that $\mathsf{E}\Gamma(T) \le \alpha + o(1)$. First, suppose that $a$ is known. Define

$$(23) \qquad\qquad \widehat{Q}_a(t) = \frac{(1-a)t}{\widehat{G}(t)}$$

to be the estimator of $Q$ when $a$ is known.

THEOREM 5.1.  *Assume that $a$ is known and let $\widehat{Q} = \widehat{Q}_a$. Let $t_0 = Q^{-1}(\alpha)$ and assume $G \ne U$. Then*

$$\sqrt{m}(T_{\text{PI}} - t_0) \rightsquigarrow N(0, K_{Q^{-1}}(t_0, t_0)),$$

$$\sqrt{m}(Q(T_{\text{PI}}) - \alpha) \rightsquigarrow N(0, (Q'(t_0))^2 K_{Q^{-1}}(t_0, t_0)),$$

*and*

$$\mathsf{E}\Gamma(T_{\text{PI}}) = \alpha + o(1).$$

PROOF. The first two statements follow from Theorem 4.4 and the delta method.

For the last claim, let $0 < \delta < t_0$, write $T = T_{PI}$ and note that

$$
\begin{aligned}
|\Gamma_m(T) - \alpha| &\leq |\Gamma_m(T) - Q(T)| + |Q(T) - \alpha| \\
&\leq \sup_t |\Gamma_m(t) - Q(t)| \mathbb{1}\{T < \delta\} \\
&\quad + \sup_t |\Gamma_m(t) - Q(t)| \mathbb{1}\{T \geq \delta\} + |Q(T) - \alpha| \\
&\leq \mathbb{1}\{T < \delta\} + \sup_{t \geq \delta} |\Gamma_m(t) - Q(t)| + |Q(T) - \alpha| \\
&= \mathbb{1}\{T < \delta\} + \frac{1}{\sqrt{m}} \sup_{t \geq \delta} |\sqrt{m}(\Gamma_m(t) - Q(t))| + |Q(T) - \alpha| \\
&= O_P(m^{-1/2}).
\end{aligned}
$$

Because $0 \leq \Gamma_m \leq 1$, the sequence is uniformly integrable, and the result follows. $\square$

Next, we consider the case where $a$ is unknown and possibly nonidentifiable. In this case, as we have seen, one can still construct an estimator that is consistent for some value $a_0 < a$.

THEOREM 5.2 (Asymptotic validity of plug-in method). *Assume that $G$ is concave. Let $T = t(\hat{a}, \widehat{G})$ be a plug-in threshold where $\widehat{G}$ is the empirical CDF or the LCM and $\hat{a} \xrightarrow{P} a_0$ for some $a_0 < a$. Then*

$$
\mathsf{E}\Gamma(T) \leq \alpha + o(1).
$$

PROOF. First note that the concavity of $G$ implies that $Q(t) = (1 - a)t/G(t)$ is increasing. Let $\delta = (a - a_0)/(1 - a)$ so that $(1 - a_0)/(1 - a) = 1 + \delta$. Then

$$
\begin{aligned}
\widehat{Q}(t) &= \frac{(1 - \hat{a})t}{\widehat{G}(t)} = \frac{1 - \hat{a}}{1 - a_0}(1 + \delta)\widehat{Q}_a(t) \\
&= (1 + o_P(1))(1 + \delta)\widehat{Q}_a(t),
\end{aligned}
$$

where $\widehat{Q}_a$ is defined in (23). Hence

$$
\begin{aligned}
T &= \widehat{Q}^{-1}(\alpha) = \widehat{Q}_a^{-1}\left(\frac{\alpha}{1 + \delta} + o_P(1)\right) \\
&\leq \widehat{Q}_a^{-1}(\alpha + o_P(1)) = \widehat{Q}_a^{-1}(\alpha) + o_P(1).
\end{aligned}
$$

Because $\widehat{Q}^{-1} \stackrel{\text{a.s.}}{\to} Q_{a_0}^{-1}$ and because $Q_{a_0}^{-1}(\alpha) \leq Q_{\underline{a}}^{-1}(\alpha)$, the result follows from the argument used in the proof of the previous theorem using $Q_{a_0}$ in place of $Q_a$. □

Recall that the oracle procedure is defined by $T_O(P^m) = Q^{-1}(\alpha)$. This procedure has the smallest FNR for all procedures that attain FDR $\leq \alpha$ up to sets of exponentially small probability [cf. Genovese and Wasserman (2002), page 506]. In the nonidentifiable case, no data-based method can distinguish $a$ and $\underline{a}$, so the performance of this oracle cannot be attained. We thus define the achievable oracle procedure $T_{AO}$ to be analogous to $T_O$ with $(1 - \underline{a})t/G(t)$ replacing $Q$. The plug-in procedure that uses the estimator $\widehat{a}$ described in Theorem 3.2 asymptotically attains the performance of $T_{AO}$ in the sense that $\mathsf{E}\Gamma(T_{PI}) = \alpha + o(1)$ and $\mathsf{E}\Xi(T_{PI}) = \mathsf{E}\Xi(T_{AO}) + o(1)$.

## 6. Confidence envelopes for FDP.
Because the distribution of the FDP need not be concentrated around its expected value, controlling the FDR does not necessarily offer high confidence that the FDP will be small. As an alternative, we develop methods in this section for making inferences about the FDP process.

A $1 - \alpha$ confidence envelope for the FDP process is a random function $\overline{\Gamma}$ on $[0,1]$ such that

$$\mathsf{P}\{\Gamma(t) \leq \overline{\Gamma}(t) \text{ for all } t\} \geq 1 - \alpha.$$

In this section we give two methods for constructing such a $\overline{\Gamma}$, one asymptotic, one exact in finite samples. See also Havránek and Chytil (1983), Hommel and Hoffman (1987) and Halperin, Lan and Hamdy (1988).

Besides being informative in its own right, a confidence envelope can be used to construct thresholds that control quantiles of the FDP distribution. We call $T$ a $1 - \alpha$ *confidence threshold* if there exists a statistic $Z$ such that

$$\mathsf{P}\{\Gamma(T) \leq Z\} \geq 1 - \alpha.$$

We consider two cases. In the first, called *rate ceiling confidence thresholds*, we take $Z$ to be a prespecified constant $c$ (the ceiling). The thresholds we develop here are derived from a confidence envelope $\overline{\Gamma}$ as the maximal threshold such that $\overline{\Gamma} \leq c$. In the second, called *minimum rate confidence thresholds*, the threshold is derived from $\overline{\Gamma}$ by $T = \arg\min_t \overline{\Gamma}(t)$ and $Z = \overline{\Gamma}(T)$.

When $a$ is known, it is possible to construct an asymptotic rate ceiling confidence threshold directly.

THEOREM 6.1. *Let* $t_c = Q^{-1}(c)$ *and let* $K_\Omega(s,t)$ *be the covariance kernel defined in* (22). *Assume that* $F \neq U$. *Define*

$$t_{c,m} \equiv t_{c,m}(\alpha) = t_c - \frac{z_\alpha}{\sqrt{m}} \frac{\sqrt{K_\Omega(t_c, t_c)}}{1 - a - cg(t_c)}.$$

*Then*

$$\mathsf{P}\{\Gamma(t_{c,m}) \le c\} = 1 - \alpha + O(m^{-1/2}).$$

PROOF. We have

$$\mathsf{P}\{\Gamma(t_{c,m}) \le c\} = \mathsf{P}\{\Omega_c(t_{c,m}) - \mu(t_{c,m}) \le -\mu(t_{c,m})\}$$

$$= \mathsf{P}\left\{\sqrt{m}\frac{\Omega_c(t_c)}{\sqrt{K_\Omega(t_c, t_c)}} \le -\frac{\sqrt{m}\mu(t_{c,m})}{\sqrt{K_\Omega(t_c, t_c)}}\right\} + o(1),$$

from Lemma 6.1. It suffices, in light of Theorem 4.1 and Lemma 6.1, to show that

$$-\sqrt{m}\frac{\mu(t_{c,m})}{\sqrt{K_\Omega(t_c, t_c)}} \to z_\alpha.$$

Now $\mu(t_c) = (1-a)t_c - cG(t_c) = 0$ since $Q(t_c) = c$. Hence

$$\mu(t) = (t - t_c)\mu'(t_c) + o(|t - t_c|)$$

$$= (t - t_c)(1 - a - cg(t_c)) + o(|t - t_c|).$$

Hence

$$\mu(t_{c,m}) = (t_{c,m} - t_c)(1 - a - cg(t_c)) + o(m^{-1/2}).$$

The result follows from the definition of $t_{c,m}$. □

LEMMA 6.1. *Let* $t_c = Q^{-1}(c)$, *and assume* $0 < t_c < 1$. *If* $t_{c,m} - t_c = O(m^{-1/2})$, $\Omega_c(t_{c,m}) - \mu(t_{c,m}) = \Omega_c(t_c) + o_P(m^{-1/2})$. *Thus, if* $u_m = vm^{-1/2} + o(m^{-1/2})$ *for some* $v$,

$$\mathsf{P}\{\Omega_c(t_{c,m}) \le \mu(t_{c,m}) + u_m\} - \mathsf{P}\{\Omega_c(t_c) \le u_m\} = o(1).$$

PROOF. Note that $\mu(t_c) = (1-a)t_c - cG(t_c) = 0$ and that

$$|\Omega_c(t_{c,m}) - \Omega_c(t_c)| \le \max\{c, 1-c\}m^{-1}\sum_i |\mathbb{1}\{P_i \le t_{c,m}\} - \mathbb{1}\{P_i \le t_c\}|$$

$$\le |\widehat{G}(t_{c,m}) - \widehat{G}(t_c)|,$$

which is Binomial$(m, |G(t_{c,m}) - G(t_c)|)/m$ and has variance of order $m^{-3/2}$. Hence

$$\Omega_c(t_{c,m}) - \mu(t_{c,m}) - \Omega_c(t_c)$$

$$= \Omega_c(t_{c,m}) - \mu(t_{c,m}) - \Omega_c(t_c) - (\mu(t_{c,m}) - \mu(t_c)) + (\mu(t_{c,m}) - \mu(t_c))$$

$$= O_P\left(\frac{1}{m^{3/4}}\right) - \mu(t_c)$$

$$= O_P\left(\frac{1}{m^{3/4}}\right) = o_P\left(\frac{1}{\sqrt{m}}\right).$$

The second claim is immediate.  □

However, when $a$ is unknown, there is a problem. When plugging in a consistent estimator of $a$ that converges at a sub-$\sqrt{m}$ rate, the error in $\hat{a}$ is of larger order than $t_c - t_{c,m}$. Using an estimator, such as Storey's estimator, which converges at a $1/\sqrt{m}$ rate but is asymptotically biased, causes overcoverage because the asymptotic bias dominates. Interestingly, as demonstrated in the next section, it is possible to ameliorate the bias problem, but not the rate problem, with appropriate conditions. Thus, a "better" estimator of $a$ need not lead to a valid confidence threshold.

6.1. *Asymptotic confidence envelope.* In this section, we show how to obtain an asymptotic confidence envelope for $\Gamma$, centered at $\widehat{Q}$. Throughout this section we use $\widehat{G}$ based on the empirical distribution function, *not* the LCM.

For reasons explained in the last section, we use Storey's estimator rather than the consistent estimators of $a$ described earlier. That is, let $\hat{a}_0 = (\widehat{G}(t_0) - t_0)/(1 - t_0)$ be Storey's estimator for a fixed $t_0 \in (0, 1)$. Then

$$\widehat{Q}(t) = \frac{(1 - \hat{a}_0)t}{\widehat{G}(t)} = \frac{1 - \widehat{G}(t_0)}{1 - t_0} \frac{t}{\widehat{G}(t)}.$$

To make the asymptotic bias in Storey's estimator negligible, we make the additional assumption that $F$ depends on a further parameter $\nu = \nu(m)$ in such a way that

$$(24) \qquad\qquad F_\nu(t) \geq 1 - e^{-\nu c(t)}$$

for some $c(t) > 0$, for all $0 < t < 1$. The marginal distribution of $P_i$ becomes

$$G_m = (1 - a)U + aF_{\nu(m)}.$$

This assumption will hold in a variety of settings such as the following:

1. The $p$-values $P_i$ are computed from some test statistics $Z_i$ that involve a common sample size $n$, where the tests all satisfy the standard large deviation principle [van der Vaart (1998), page 209]. In this case $\nu = n$.
2. As in the previous case except that each test has a sample size $n_i$ drawn from some common distribution.
3. Each test is based on measurements from a counting process (such as an astronomical image) where $\nu$ represents exposure time.

Under these assumptions, we have the following theorem.

THEOREM 6.2. *Let $t_m$ be such that $t_m \to 0$ and $m t_m / (\log m)^4 \to \infty$. Let $w_{\alpha/2}$ denote the upper $\alpha/2$ quantile of $\max_{0 \leq t \leq 1} \mathbb{B}(t)/\sqrt{t}$, where $\mathbb{B}(t)$ denotes a standard Brownian bridge. Let*

$$\Delta_m = \max\left\{ 2(1 - \hat{a}_0) w_{\alpha/2}, \ \frac{\sqrt{2}}{1 - t_0} \sqrt{\log\left(\frac{4}{\alpha}\right)} \right\}. \tag{25}$$

*Define*

$$\overline{\Gamma}(t) = \min\left\{ \widehat{Q}(t) + \frac{\Delta_m \sqrt{t}}{\sqrt{m}\widehat{G}(t)}, 1 \right\}. \tag{26}$$

*Assume that*

$$\frac{\nu(m)}{\log m} \to \infty \tag{27}$$

*as $m \to \infty$. Then*

$$\liminf_{m \to \infty} \mathsf{P}\{\Gamma(t) \leq \overline{\Gamma}(t) \text{ for all } t \geq t_m\} \geq 1 - \alpha. \tag{28}$$

PROOF. Let

$$N(t) = \frac{M_{1|0}(t)}{m} = \frac{1}{m} \sum_{i=1}^{m} (1 - H_i) \mathbb{1}\{P_i \leq t\}.$$

Note that $\mathsf{E}(N(t)) = (1 - a)t$ and $\text{Cov}(N(t), N(s)) = (1 - a)^2(s \wedge t - st)$. By Donsker's theorem, $\sqrt{m}(N(t) - (1 - a)t) \rightsquigarrow (1 - a)\mathbb{B}(t)$, where $\mathbb{B}(t)$ is a standard Brownian bridge. By the Hungarian embedding, there exists a sequence of standard Brownian bridges $\mathbb{B}_m(t)$ such that

$$N(t) = (1 - a)t + \frac{(1 - a)\mathbb{B}_m(t)}{\sqrt{m}} + R_m(t),$$

where

$$R_m \equiv \sup_t |R_m(t)| = O\left(\frac{(\log m)^2}{m}\right) \qquad \text{a.s.}$$

Let

$$V(t) = (1 - \hat{a}_0)t + \frac{\Delta_m \sqrt{t}}{\sqrt{m}}. \tag{29}$$

Now,

$$\mathsf{P}\{N(t) > V(t) \text{ for some } t \geq t_m\}$$

$$= \mathsf{P}\left\{ (1 - a)t + \frac{(1 - a)\mathbb{B}_m(t)}{\sqrt{m}} + R_m(t) \right.$$

$$(30) \quad \begin{aligned} &> (1 - \hat{a}_0)t + \frac{\Delta_m \sqrt{t}}{\sqrt{m}} \text{ for some } t \geq t_m \Big\} \\ &= \mathsf{P}\Big\{ \max_{t \geq t_m} \Big( \sqrt{m}(\hat{a}_0 - a)\sqrt{t} + (1 - a)\frac{\mathbb{B}_m(t)}{\sqrt{t}} + \frac{\sqrt{m}R_m}{\sqrt{t}} \Big) > \Delta_m \Big\} \\ &\leq \mathsf{P}\Big\{ \max_{t \geq t_m} (\sqrt{m}|\hat{a}_0 - a|\sqrt{t}) > \frac{\Delta_m}{2} \Big\} + \mathsf{P}\Big\{ (1 - a)\max_{t \geq t_m} \frac{\mathbb{B}_m(t)}{\sqrt{t}} > \frac{\Delta_m}{2} \Big\} \\ &\quad + O\Big( \frac{(\log m)^2}{\sqrt{t_m}\sqrt{m}} \Big). \end{aligned}$$

The last term is $o(1)$ since $mt_m/(\log m)^4 \to \infty$.

Let

$$a_0 = \frac{G(t_0) - t_0}{1 - t_0} = a \frac{F_{\nu(m)}(t_0) - t_0}{1 - t_0}.$$

Then

$$a - a_0 = a \frac{1 - F_{\nu(m)}(t_0)}{1 - t_0} \leq \frac{e^{-\nu(m)c(t_0)}}{1 - t_0}.$$

By assumption, we can write

$$\nu(m) = \frac{s_m \log m}{c(t_0)}$$

for some $s_m \to \infty$. Hence $a - a_0 = O(m^{-s_m})$. In particular, $a - a_0 = o(\frac{1}{\sqrt{m}})$. Hence

$$\sqrt{m}|\hat{a}_0 - a| \leq \sqrt{m}|\hat{a}_0 - a_0| + \sqrt{m}|a_0 - a| = \sqrt{m}|\hat{a}_0 - a_0| + o(1).$$

Thus

$$\begin{aligned} &\mathsf{P}\Big\{ \max_{t \geq t_m} (\sqrt{m}|\hat{a}_0 - a|\sqrt{t}) > \frac{\Delta_m}{2} \Big\} \\ &= \mathsf{P}\Big\{ \sqrt{m}|\hat{a}_0 - a| > \frac{\Delta_m}{2} \Big\} \\ &= \mathsf{P}\Big\{ \sqrt{m}|\hat{a}_0 - a_0| > \frac{\Delta_m}{2} \Big\} + o(1) \\ &= \mathsf{P}\Big\{ \frac{\sqrt{m}|\hat{G}(t_0) - G_m(t_0)|}{1 - t_0} > \frac{\Delta_m}{2} \Big\} + o(1) \\ &= \mathsf{P}\Big\{ |\hat{G}(t_0) - G_m(t_0)| > \frac{\Delta_m(1 - t_0)}{2\sqrt{m}} \Big\} + o(1) \\ &\leq 2\exp\Big\{ -m\frac{2\Delta_m^2}{4}\frac{(1 - t_0)^2}{m} \Big\} + o(1) \\ (31) \quad &\leq \frac{\alpha}{2} + o(1) \end{aligned}$$

by the DKW inequality and the definition of $\Delta_m$.

Fix $\varepsilon > 0$. Since $\hat{a}_0 \overset{\text{a.s.}}{\to} a_0$, we have, almost surely for all large $m$, that

$$\frac{\Delta_m}{2(1-a)} \geq \frac{2(1-\hat{a}_0)w_{\alpha/2}}{2(1-\hat{a})}$$

$$= \frac{1-\hat{a}_0}{1-a}w_{\alpha/2} = \frac{1-\hat{a}_0}{1-a_0}(1+o(1))w_{\alpha/2} \geq w_{\alpha/2} - \varepsilon.$$

Let $\mathbb{W}_m(t) = \mathbb{B}_m(t)/\sqrt{t}$. Then for all large $m$

$$\mathsf{P}\left\{(1-a)\max_{t \geq t_m}\mathbb{W}_m(t) > \frac{\Delta_m}{2}\right\}$$

$$= \mathsf{P}\left\{\max_{t \geq t_m}\mathbb{W}_m(t) > \frac{\Delta_m}{2(1-a)}\right\}$$

$$\leq \mathsf{P}\left\{\max_{t \geq t_m}\mathbb{W}_m(t) > w_{\alpha/2} - \varepsilon\right\}$$

$$\leq \mathsf{P}\left\{\max_{0 \leq t \leq 1}\mathbb{W}_m(t) > w_{\alpha/2} - \varepsilon\right\}$$

$$= \mathsf{P}\left\{\max_{0 \leq t \leq 1}\mathbb{W}_m(t) > w_{\alpha/2}\right\}$$

$$+ \mathsf{P}\left\{w_{\alpha/2} - \varepsilon < \max_{0 \leq t \leq 1}\mathbb{W}_m(t) \leq w_{\alpha/2}\right\}$$

$$= \frac{\alpha}{2} + \mathsf{P}\left\{w_{\alpha/2} - \varepsilon < \max_{0 \leq t \leq 1}\mathbb{W}_m(t) \leq w_{\alpha/2}\right\}.$$

Since $\varepsilon$ is arbitrary, this implies that

$$(32) \qquad \limsup_{m \to \infty} \mathsf{P}\left((1-a)\max_{t \geq t_m}\mathbb{W}_m(t) > \frac{\Delta_m}{2}\right) \leq \frac{\alpha}{2}.$$

From (31), (32) and (30) we conclude that

$$\limsup_{m \to \infty} \mathsf{P}(N(t) > V(t) \text{ for some } t \geq t_m) \leq \alpha.$$

Notice that $\Gamma(t) = N(t)/\hat{G}(t)$. Hence $N(t) \leq V(t)$ implies that

$$\Gamma(t) \leq \frac{V(t)}{\hat{G}(t)} = \overline{\Gamma}(t).$$

The conclusion follows. $\quad\square$

Both types of confidence thresholds can now be defined from $\overline{\Gamma}$. For example, pick a ceiling $0 < c < 1$ and define $T_c = \sup\{t \geq t_m : \overline{\Gamma}(t) \leq c\}$, where $T_c$ is defined to be 0 if no such $t$ exists. The proof of the following is then immediate from the previous theorem.

COROLLARY 6.1. *$T_c$ is an asymptotic rate ceiling confidence threshold with ceiling $c$.*

It is also worth noting that we can construct a confidence envelope for the number of false discoveries process $M_{1|0}(t)$.

COROLLARY 6.2. *With $t_m$ as in the above theorem and $V(t)$ defined as in* (29),

$$\liminf_{m\to\infty} \mathsf{P}\{M_{1|0}(t) \le mV(t) \text{ for } t \ge t_m\} \ge 1 - \alpha. \tag{33}$$

6.2. *Exact confidence envelope.* In this section we will construct confidence thresholds that are valid for finite samples.

Let $0 < \alpha < 1$. Given $V_1, \dots, V_k$, let $\varphi_k(v_1, \dots, v_k)$ be a nonrandomized level $\alpha$ test of the null hypothesis that $V_1, \dots, V_k$ are drawn i.i.d. from a Uniform$(0,1)$ distribution. Define $p_0^m(h^m) = (p_i : h_i = 0, 1 \le i \le m)$ and $m_0(h^m) = \sum_{i=1}^m (1 - h_i)$ and $\mathcal{U}_\alpha(p^m) = \{h^m \in \{0,1\}^m : \varphi_{m_0(h^m)}(p_0^m(h^m)) = 0\}$. Note that as defined, $\mathcal{U}_\alpha$ always contains the vector $(1, 1, \dots, 1)$.

Let

$$\mathcal{G}_\alpha(p^m) = \{\Gamma(\cdot, h^m, p^m) : h^m \in \mathcal{U}_\alpha(p^m)\}, \tag{34}$$

$$\mathcal{M}_\alpha(p^m) = \{m_0(h^m) : h^m \in \mathcal{U}_\alpha(p^m)\}. \tag{35}$$

Then we have the following theorem, which follows from standard results on inverting hypothesis tests to construct confidence sets.

THEOREM 6.3. *For all $0 < a < 1, F$, and positive integers $m$,*

$$\mathsf{P}_{a,F}\{H^m \in \mathcal{U}_\alpha(P^m)\} \ge 1 - \alpha, \tag{36}$$

$$\mathsf{P}_{a,F}\{M_0 \in \mathcal{M}_\alpha(P^m)\} \ge 1 - \alpha, \tag{37}$$

$$\mathsf{P}_{a,F}\{\Gamma(\cdot, H^m, P^m) \in \mathcal{G}_\alpha\} \ge 1 - \alpha, \tag{38}$$

$$\mathsf{P}_{a,F}\{\Gamma(T_c) \le c\} \ge 1 - \alpha, \tag{39}$$

*where*

$$T_c = \sup\{t : \Gamma(t; h^m, P^m) \le c \text{ and } h^m \in \mathcal{U}_\alpha(P^m)\}. \tag{40}$$

*In particular,*

$$\overline{\Gamma}(t) = \sup\{\Gamma(t) : \Gamma \in \mathcal{G}_\alpha(P^m)\} \tag{41}$$

*is a $1 - \alpha$ confidence envelope for $\Gamma$, and $T_c$ is a $1 - \alpha$ rate ceiling confidence threshold with ceiling $c$. In fact, $\inf_{a,F} \mathsf{P}_{a,F}\{\Gamma(t) \le \overline{\Gamma}(t), \text{ for all } t\} \ge 1 - \alpha$.*

REMARK 6.1.   If there is some substantive reason to bound $M_0$ from below, then $\mathcal{G}_\alpha$ will have a nontrivial lower bound as well. In general, because $\mathcal{U}_\alpha$ always contains $(1, 1, \ldots, 1)$, the pointwise infimum of functions in $\mathcal{G}_\alpha$ will be zero.

REMARK 6.2.   At first glance, computation of $\mathcal{U}_\alpha$ would appear to require an exponential-time algorithm. However, for broad classes of tests, including the Kolmogorov–Smirnov test, it is possible to construct $\mathcal{U}_\alpha$ in polynomial time.

REMARK 6.3.   The choice of test can be important for obtaining a good confidence envelope. A full analysis of this choice is beyond the scope of this paper; we will present such an analysis in a forthcoming paper. In the examples below, we use the test derived from the second-order statistic of a subset of $p$-values.

REMARK 6.4.   A similar construct yields a confidence envelope on the process $M_{1|0}(t)$.

### 6.3. *Examples.*

EXAMPLE 1.   We begin with a re-analysis of Example 3.2 from BH (1995). BH give the following 15 $p$-values

| 0.0001 | 0.0004 | 0.0019 | 0.0095 | 0.0201 | 0.0278 | 0.0298 | 0.0344 |
| 0.0459 | 0.3240 | 0.4262 | 0.5719 | 0.6528 | 0.7590 | 1      |        |

and at a 0.05 level Bonferroni rejects the first three null hypotheses and the BH method rejects the first four.

Because $m$ is small, we construct only the exact confidence envelope for this example. Figure 1 shows the upper 95% confidence envelope on the FDP for these data using the second-order statistic of any subset as a test statistic for the exact procedure. Notice first that the confidence envelope never drops below 0.05. Second, while the BH threshold $T = P_{(4)} = 0.0095$ guarantees an FDR $\leq 0.05$, we can claim that $\mathsf{P}\{\Gamma(P_{(4)}) > 0.25\} \leq 0.05$, but this is also true for the larger threshold $P_{(11)}^- = 0.4262^-$, which will have higher power. This difference occurs because the envelope takes large values at small thresholds. The result could be quite different with another choice of test statistic. The minimum rate 95% confidence threshold has $T = 0.324$ and $Z = \overline{\Gamma}(T) = 0.111$.

EXAMPLE 2. We present a simple, synthetic example, where $m = 1000$, $a = 0.25$, and the test-statistic is from a Normal$(\theta, 1)$ one-sided test with $H_0 : \theta = 0$ and $H_1 : \theta = 3$.

Figure 2 compares the true FDP sample path with the 95% confidence envelopes derived from the exact and asymptotic methods. For small values of the threshold the exact envelope almost matches the truth, but for larger values it becomes more conservative. The asymptotic envelope remains above but generally close to the truth. The asymptotic and exact envelopes cross at an FDP level of about 0.05. The rate ceiling confidence thresholds with ceiling 0.05 and level 0.05 are 0.00062 for the asymptotic and 0.00046 for the exact. The minimum rate confidence threshold for the exact procedure has $T = 0.00039$ and $Z = 0.011$.

## APPENDIX

**Algorithm for finding $\widehat{F}_m$.** Here we restrict our attention to the case in which we take $\hat{F}$ as piecewise constant on the same grid as $\mathbb{G}$. When $F$ is concave, the algorithm works in the same way with the sharper piecewise linear approximation.

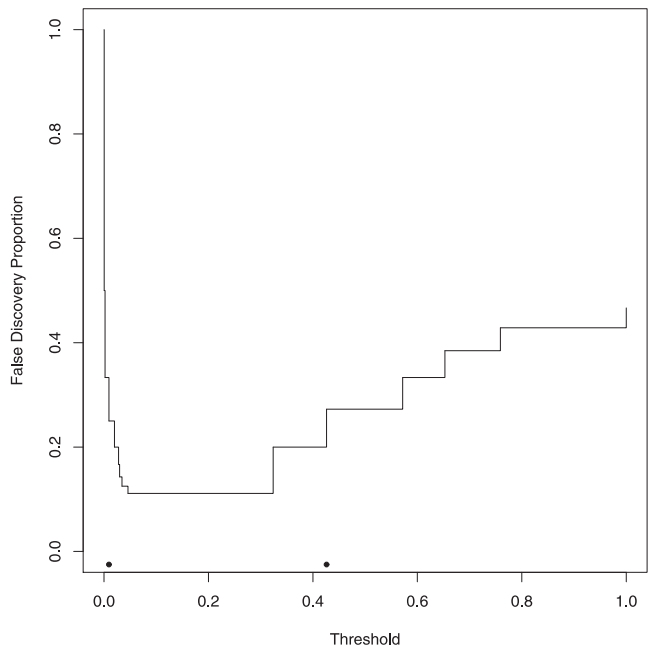

FIG. 1. *Plot of $\overline{\Gamma}(t)$ versus $t$ for Example 1, where $\overline{\Gamma}$ is derived from the exact method of Section 6.2. The leftmost dot on the horizontal axis is the BH threshold; the rightmost dot is a confidence threshold with the same ceiling.*

STEP 0. Begin by constructing an initial estimate of $F$ that is a CDF. For example, we can define $H$ to be the piecewise constant function that takes the following values on the $P_i$'s:

$$H(P_{(i)}) = \max_{j \le i} \frac{\widehat{G}(P_{(j)}) - (1 - \hat{a})P_{(j)}}{\hat{a}}.$$

STEP 1. Identify the segment with the biggest absolute difference between $\widehat{G}$ and $(1 - \hat{a})U + \hat{a}H$.

STEP 2. Determine how far and in what direction (up or down) this segment can be moved while keeping $H$ a CDF and minimizing $\|\widehat{G} - (1 - \hat{a})U + \hat{a}H\|_\infty$.

STEP 3. If the segment can be moved, move it and go to Step 1. Else go to Step 4.

STEP 4. If no segment can be moved to reduce $\|\widehat{G} - (1 - \hat{a})U + \hat{a}H\|_\infty$, STOP.

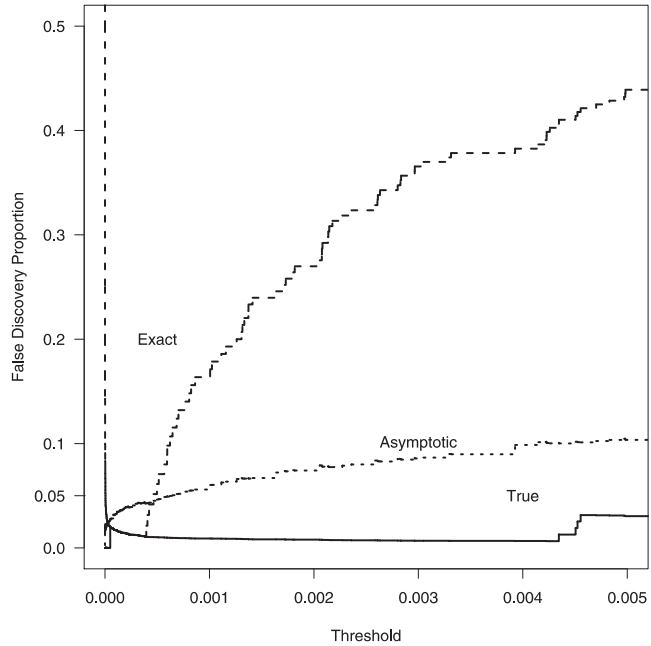

FIG. 2. *Plot of the true $\Gamma$ sample paths and $\overline{\Gamma}$ for the exact (cf. Section 6.2) and asymptotic (cf. Section 6.1) methods for the data in Example 2. The envelopes are shown here only for small thresholds. The truth (solid) is the lowest curve over the entire domain. The exact envelope (dashed) begins near 1, dips toward the truth and then rises sharply. The asymptotic envelope (dotted) is the other curve.*

If the current segment is part of a contiguous block of segments where one segment in the block can be moved to reduce $\|\widehat{G} - (1 - \hat{a})U + \hat{a}H\|_\infty$, move the segment at the end of the contiguous block of segments that provides the greatest reduction in $\|\widehat{G} - (1 - \hat{a})U + \hat{a}H\|_\infty$. Go to Step 1.

**Acknowledgments.** The authors are grateful to the referees for providing many helpful suggestions.

Department of Statistics
Carnegie Mellon University
Pittsburgh, Pennsylvania 15213
USA
e-mail: larry@stat.cmu.edu